\documentclass[a4, 10pt]{amsart}
\usepackage{amssymb}
\usepackage{amstext}
\usepackage{amsmath}
\usepackage{amscd}
\usepackage{amsthm}
\usepackage{latexsym}
\usepackage{amsfonts}
\usepackage{graphicx}
\usepackage{enumerate}
\usepackage{color}

\theoremstyle{plain}
\newtheorem{thm}{Theorem}[section]
\newtheorem{prop}[thm]{Proposition}
\newtheorem{lemma}[thm]{Lemma}
\newtheorem{cor}[thm]{Corollary}

\theoremstyle{definition}

\newtheorem{ex}[thm]{Example}

\theoremstyle{remark}
\newtheorem*{pf}{{\sl Proof}}
\newtheorem*{pf1}{{\sl Proof of Theorem \ref{tf}}}

\numberwithin{equation}{section}

\def\Hom{\mathrm{Hom}}

\def\Ext{\mathrm{Ext}}
\def\Tor{\mathrm{Tor}}

\def\mod{\mathrm{mod}}
\def\Coker{\mathrm{Coker}}

\def\Im{\mathrm{Im}}
\def\tr{\mathrm{Tr}}

\def\m{\mathfrak m}

\def\p{\mathfrak p}
\def\q{\mathfrak q}

\def\Z{\Bbb Z}

\def\depth{\mathrm{depth\,}}

\def\id{\mathrm{id}}

\def\height{\mathrm{ht}}

\def\Spec{\mathrm{Spec}}

\def\mod{\mathrm{mod}}

\begin{document}

\title{Remarks on torsionfreeness and its applications}

\author{Tokuji Araya}
\address{Department of Applied Science, Faculty of Science, Okayama University of Science, Ridaicho, Kitaku, Okayama 700-0005, Japan.}
\email{araya@das.ous.ac.jp}

\author{Kei-ichiro Iima}
\address{Department of Liberal Studies, National Institute of Technology, Nara College, 22 Yata-cho, Yamatokoriyama, Nara 639-1080, Japan.}
\email{iima@libe.nara-k.ac.jp}

\thanks{The first author was partially supported by JSPS Grant-in-Aid for Scientific Research (C) 26400056}
\maketitle


\begin{abstract}
In this article, we shall characterize torsionfreeness of modules with respect to a semidualizing module in terms of the Serre's condition $(S_n)$.
As its applications, we give a characterization of Cohen-Macaulay rings $R$ such that $R_\p$ is Gorenstein for all prime ideals $\p$ of height less than $n$, and we will give a partial answer of Tachikawa conjecture and Auslander-Reiten conjecture.
\end{abstract}



\section{Introduction}

There are many researches about maximal Cohen-Macaulay modules over Cohen-Macaulay rings from the viewpoint of Cohen-Macaulay representation theory.
For instance, one can easily see that all modules over zero dimensional rings are maximal Cohen-Macaulay. 
It is well known results that a torsionless (resp. reflexive) module over one dimensional domains (resp. two dimensional normal domains) is coincide with a maximal Cohen-Macaulay module (c.f. \cite{BH, Y}).
On the other hands, one can check that a second syzygy module of the residue class field of three dimensional normal Cohen-Macaulay local domain is not maximal Cohen-Macaulay but reflexive. 
Auslander and Bridger introduce a notion of $n$-$torsionfree$ as generalization of reflexive \cite{AB}. 
Evans and Griffith give a characterization of $n$-torsionfree modules \cite{EG}. 

A $semidualizing$ $module$ has been defined by Foxby \cite{F}, Vasconcelos \cite{V} and Golod \cite{G}.
A free module of rank one and a canonical module over Cohen-Macaulay rings are typical examples of semidualizing modules.

The notion of $n$-$C$-$torsionfree$ has been introduced by Takahashi \cite{T}.
In this article, we study an $n$-$C$-torsionfreeness of modules with respect to a semidualizing module in terms of the Serre's condition $(S_n)$.
Recently, Dibaei and Sadeghi \cite{DS} give a similar property, independently.

\begin{thm}\label{tf}
Let $R$ be a commutative noetherian ring and let $M$ be an $R$-module. Assume that $R$ satisfies the Serre's condition $(S_n)$.
Furthermore, we assume that either $R$ satisfies $(G_{n-1}^C)$ or $M_\p$ is totally $C_\p$-reflexive for each prime ideal $\p$ of height less than $n$. 
Then the following conditions are equivalent for $M$:
\begin{enumerate}[\rm (1)]
\item $M$ is $n$-$C$-torsionfree,
\item $M$ has an $n$-$C$-universal pushforward,
\item $M$ is $n$-$C$-syzygy,
\item $M$ satisfies $(\widetilde{S_n})$,
\item $M$ satisfies $(S_n)$.
\end{enumerate}
\end{thm}

In Section $2$, we prepare some notions and some properties to explain Theorem \ref{tf}. 
In Section $3$, we prove Theorem \ref{tf}. 
In Section $4$, we give applications of Theorem \ref{tf}.
The first application is to characterize Cohen-Macaulay rings $R$ such that $R_\p$ is Gorenstein for all prime ideals $\p$ of height less than $n$. 
Tachikawa conjecture \cite{Tac} and Auslander-Reiten conjecture \cite{AR} are long-standing conjectures which originated in Nakayama conjecture \cite{Nak}.
We attempt to solve Tachikawa conjecture and Auslander-Reiten conjecture as a second application.
In Section $5$, we give some examples of non-trivial semidualizing modules.


\section{Preliminaries}

Throughout the rest of this article, let $R$ be a commutative noetherian ring.
All modules are assumed to be finitely generated.
In this section, we give some notions and properties with respect to a semidualizing module.

An $R$-module $C$ is called {\it semidualizing} if the homothety map $R\to \Hom_R(C,C)$ is an isomorphism and if $\Ext^{>0}_R(C,C)=0$.
From now on, we fix a semidualizing module $C$ and put $(-)^{\dag}=\Hom_R(-,C)$.

Let $\cdots \to P_1 \overset{\partial}{\to} P_0 \to M \to 0$ be a projective resolution of an $R$-module $M$.
We define a {\it $C$-transpose} module $\tr_CM$ of $M$ the cokernel of $P_0^{\dag}\overset{\partial^{\dag}}{\to} P_1^{\dag}$. 
We remark that $\tr_CM$ is uniquely determined up to direct summands of finite direct sums of copy of $C$. 
Note that if $C$ is isomorphic to $R$ then $C$-transpose coincides with ordinary (Auslander) transpose.
An $R$-module $M$ is called {\it $n$-$C$-torsionfree} if $\Ext^i_R(\tr_CM,C)=0$ for all $1 \leq i \leq n$.

We denote by $\lambda_M$ the natural map $M\to M^{\dag\dag}$.
$n$-$C$-torsionfreeness has following properties similar to ordinary $n$-torsionfreeness \cite{AB}.
One can show this by diagram chasing (c.f. \cite{T}).

\begin{prop}\label{2.1}
Let $M$ be an $R$-module.
\begin{enumerate}[\rm (1)]
\item $M$ is $1$-$C$-torsionfree if and only if $\lambda_M$ is a monomorphism,
\item $M$ is $2$-$C$-torsionfree if and only if $\lambda_M$ is an isomorphism,
\item Let $n\geq 3$. $M$ is $n$-$C$-torsionfree if and only if $\lambda_M$ is an isomorphism and if $\Ext^i_R(M^{\dag},C)=0$ for all $1 \leq i \leq n-2$.
\end{enumerate}
\end{prop}

An $R$-module $M$ is called {\it totally $C$-reflexive} if the natural map $\lambda_M$ is an isomorphism and if ${\operatorname{Ext}}^{>0}_R(M,C)=0=\Ext^{>0}_R(M^{\dag},C)$.
Note that if $C$ is isomorphic to $R$ then the notion of totally $C$-reflexive coincides with the notion of ordinary totally reflexive.
Each condition (2), (3) and (4) in the next lemma gives a characterization of totally $C$-reflexive modules.
They are well-known if $C=R$.
One can show this lemma similar to the case where $C=R$.
This lemma plays key role to show Theorem \ref{main}.

\begin{lemma}
The following conditions are equivalent for an $R$-module $M$:
\begin{enumerate}[\rm (1)]
\item $M$ is totally $C$-reflexive.
\item $M_\p$ is totally $C_\p$-reflexive for each prime ideal $\p$.
\item There exists an exact sequence $\cdots \to P^{-1} \to P^{0} \overset{\partial}{\to} P^{1}_C \to P^{2}_C \to\cdots$ such that $M$ is isomorphic to $\Im \partial$, $P^i$ is projective for each $i\leq 0$, $P^i_C$ is a direct summand of finite direct sum of copies of $C$  for each $i>1$, and $C$-dual sequence $\cdots \to (P^{2}_C)^{\dag} \to (P^{1}_C)^{\dag} \to (P^{0})^{\dag} \to (P^{-1})^{\dag} \to\cdots$ is exact.
\item The syzygy module $\Omega M$ of $M$ is totally $C$-reflexive and $\Ext^1_R(M,C)=0$.
\end{enumerate}
When this is the case, we have $\depth M=\depth R$ if $R$ is local.
\end{lemma}

An $R$-module $M$ is called $n$-$C$-{\it syzygy} if there exists an exact sequence $0 \to M \to P_C^1 \to P_C^2 \to \cdots \to P_C^n$ such that each $P_C^i$ is a direct summand of finite direct sums of copy of $C$.
We set $\Omega_C^n(\mod R)$ the class of $n$-$C$-syzygy modules. 
Let $\sigma$ : $0 \to M \to P_C^1 \to P_C^2 \to \cdots \to P_C^n$ be a exact sequence whose $P_C^i$ is a direct summand of finite direct sums of copy of $C$ for each $i$.
We say $\sigma$ $n$-$C$-{\it universal pushforward} of $M$ if the $C$-dual sequence $\sigma^{\dag}$ is exact. 
Note that if $M$ has an $n$-$C$-universal pushforward, then $M$ is $n$-$C$-syzygy.

We say that an $R$-module $M$ satisfies the Serre's condition $(S_n)$ (resp. $(\widetilde{S_n})$) if $\depth M_{\p}\geq \min\{ n, \dim R_{\p}\}$ (resp. $\depth M_{\p}\geq \min\{ n, \depth R_{\p}\}$) for each prime ideal $\p$ of $R$.
We denote by $S_n(R)$ the class of modules which satisfies $(S_n)$. 
We say that $R$ satisfies the condition $(ES^C_n)$ provided each $R$-module $X$ with $\Ext^{>0}_R(X,C)=0$ satisfies $(\widetilde{S_n})$.
Moreover, we say that $R$ satisfies the condition $(G_n^C)$ if injective dimension of $C_\p$ (as an $R_\p$-module) is finite for all prime ideal $\p$ of height at most $n$.
In this case, $R_\p$ is Cohen-Macaulay local ring with canonical module $C_\p$ for all prime ideal $\p$ of height at most $n$.


\section{Proof of Theorem \ref{tf}}

In this section, we give a proof of Theorem \ref{tf}.
We can show the equivalence between $(1)$ and $(2)$ in the Theorem \ref{tf} with no assumptions.

\begin{prop}\label{3.1}
The following two conditions are equivalent for an $R$-module $M$:
\begin{enumerate}[\rm (1)]
\item $M$ is $n$-$C$-torsionfree,
\item $M$ has an $n$-$C$-universal pushforward.
\end{enumerate}
\end{prop}

\begin{pf}
$(1)\Rightarrow (2)$ We prove this by induction on $n$.
Assume $n=1$.
Let $f_1,f_2,\ldots,f_m$ be a generating system of $M^{\dag}$ and put $f=(f_1\ f_2\ \ldots\ f_m)$.
Applying the $C$-dual to an exact sequence $R^{\oplus m} \overset{f}{\longrightarrow} M^{\dag} \to 0$, we obtain an exact sequence $0 \longrightarrow M^{\dag\dag} \overset{f^{\dag}}{\longrightarrow} C^{\oplus m}$.
Since $M$ is $1$-$C$-torsionfree, $\lambda_M$ is a monomorphism.
Putting $g=f^{\dag}\lambda_M$, we see that $0\to M\overset{g}{\to}C^{\oplus m}$ is a $1$-$C$-universal pushforward of $M$.

Assume $n \geq 2$.
Since $M$ is $1$-$C$-torsionfree, $M$ has a $1$-$C$-universal pushforward $0\to M\overset{g}{\to} P_C^1$.
Putting $N=\Coker g$, we have a following commutative diagram:

$$
\begin{CD}
0 @>>> M @>>> P_C^1 @>>> N @>>> 0\\
@. @V{\lambda_M}VV @V{\lambda_{P_C^1}}VV @V{\lambda_N}VV \\
0 @>>> M^{\dag\dag} @>>> {P_C^1}^{\dag\dag} @>>> N^{\dag\dag} @>>> \Ext^1_R(M^{\dag},C) @>>> 0.\\
\end{CD}
$$

Since $\Ext^i_R(N^{\dag},C)\cong \Ext^{i+1}_R(M^{\dag},C)$ for each $i>0$, $N$ is $(n-1)$-$C$-torsionfree.
By induction assumption, $N$ has an $(n-1)$-$C$-universal pushforward $0\to N\to P_C^2\to\dots\to P_C^n$.
Then, $0\to M\to P_C^1\to P_C^2\to\dots\to P_C^n$ is an $n$-$C$-universal pushforward of $M$.

$(2)\Rightarrow (1)$ 
Since $M$ has a $1$-$C$-universal pushforward,
there exists a short exact sequence $0\to M \to P_C\to N\to 0$ such that $P_C$ is a direct summand of finite direct sums of copy of $C$ and that the $C$-dual sequence $0\to N^{\dag}\to P_C^{\dag}\to M^{\dag}\to 0$ is exact.
Then we obtain a following commutative diagram:
$$
\begin{CD}
0 @>>> M @>>> P_C @>>> N @>>> 0\\
@. @V{\lambda_M}VV @V{\lambda_{P_C}}V{\cong}V @V{\lambda_N}VV \\
0 @>>> M^{\dag\dag} @>>> P_C^{\dag\dag} @>>> N^{\dag\dag} @>>> \Ext^1_R(M^{\dag},C) @>>> 0.\\
\end{CD}
$$
Note that $\Ext^i_R(N^{\dag},C)\cong \Ext^{i+1}_R(M^{\dag},C)$ for each $i>0$.

Now, we shall show that $M$ is $n$-$C$-torsionfree by induction on $n$.
We remark that $\lambda_M$ is a monomorphism by the five lemma.
Thus $M$ is $1$-$C$-torsionfree.

Assume $n\geq 2$.
Since $N$ has an $(n-1)$-$C$-universal pushforward, $N$ is $(n-1)$-$C$-torsionfree by induction assumption.
This yield that $M$ is $n$-$C$-torsionfree.
\qed
\end{pf}

\begin{pf1}
The implication $(2)\Rightarrow (3)$ is obvious.
Since $\depth C_{\p}=\depth R_{\p}$ for all prime ideal $\p$, one can check the implication $(3)\Rightarrow (4)$ by using depth lemma.
$R$ satisfies $(S_n)$ by the assumption, we can see that $(4)$ and $(5)$ are equivalent.
Therefore, it is enough to prove the implication $(5) \Rightarrow (1)$.
We prove the implication $(5)\Rightarrow (1)$ by induction on $n$.
Assume $n=1$.
Let $\p$ be an associated prime ideal of $M$.
Since $M$ satisfies $(S_1)$, we have $\dim R_\p=0$.
If $R$ satisfies $(G_0^C)$, then $C_\p$ is a canonical $R_\p$-module.
Therefore $M_\p$ is totally $C_\p$-reflexive and we have $\Hom_R(M,C)_\p\cong\Hom_{R_\p}(M_\p,C_\p)\not=0$.
In particular, $\Hom_R(M,C)\not= 0$.

Let $f_1,f_2,\dots,f_m$ be a generating system of $\Hom(M,C)$ and put $f=^t\hspace{-4pt}(f_1,f_2,\dots,f_m):M\to C^{\oplus m}$.
Suppose that $N=\ker f$ is not zero.
Let $\q$ be an associated prime ideal of $N$.
Since $\q$ is also an associated prime ideal of $M$, we have $\dim R_\q=0$.
Noting that $C_\q$ is canonical module over $R_\q$, we have that $f_\q$ is a monomorphism.
This yields that $N_\q=0$.
This contradicts that $\q$ is an associated prime ideal of $N$.
Hence $f$ is a monomorphism.

Since $f^{\dag\dag}\lambda_M=\lambda_{C^{\oplus m}}f$ is a monomorphism, we obtain that $\lambda_M$ is a monomorphism.
This means that $M$ is $1$-$C$-torsionfree.

Assume $n \geq 2$.
Since $M$ satisfies $(S_1)$, $M$ is $1$-$C$-torsionfree.
In particular, $M$ has a $1$-$C$-universal pushforward $0 \to M \to P_C \to N \to 0$ by Proposition \ref{3.1}.
Then we get a following commutative diagram:
$$
\begin{CD}
0 @>>> M @>>> P_C @>>> N @>>> 0\\
@. @V{\lambda_M}VV @V{\lambda_{P_C}}VV @V{\lambda_N}VV \\
0 @>>> M^{\dag\dag} @>>> P_C^{\dag\dag} @>>> N^{\dag\dag} @>>> \Ext^1_R(M^{\dag},C) @>>> 0.\\
\end{CD}
$$
Note that $\Ext^i_R(N^{\dag},C)\cong \Ext^{i+1}_R(M^{\dag},C)$ for each $i>0$.
By induction assumption, it is enough to prove that $N$ satisfies $(S_{n-1})$ and that $N_\p$ is totally $C_\p$-reflexive for all prime ideal $\p$ of height less than $n-1$.

Let $\p$ be a prime ideal.
If $\dim R_{\p}\geq n$, we have $\depth M_{\p} \geq \min \{ n, \dim R_{\p}\}=n$.
Therefore we obtain $\depth N_{\p}\geq n-1$ by depth lemma.

Assume $\dim R_{\p}\leq n-1$.
If $R$ satisfies $(G_{n-1}^C)$, $R_{\p}$ is Cohen-Macaulay with canonical module $C_{\p}$.
Inequalities $\depth M_{\p}\geq \min \{ n, \dim R_{\p}\}=\dim R_{\p}=\depth R_{\p}$ gives that $M_{\p}$ is a maximal Cohen-Macaulay $R_\p$-module.
In particular, $M_\p$ is totally $C_\p$-reflexive.
Thus so are $(M_{\p})^{\dag_\p}$, $R_{\p}$ and $(N_{\p})^{\dag_\p}$.

It comes from a commutative diagram:
$$\begin{CD}
0 @>>> M_\p @>>> (P_C)_\p @>>> N_\p @>>> 0\\
@. @V{\lambda_{M_\p}}V{\cong}V @V{\lambda_{(P_C)_\p}}V{\cong}V @V{\lambda_{N_\p}}VV \\
0 @>>> (M_\p)^{\dag_\p \dag_\p} @>>> (P_C)_\p^{\dag_\p\dag_\p} @>>> (N_\p)^{\dag_\p \dag_\p} @>>> 0,\\
\end{CD}
$$
we can see that $\lambda_{N_\p}$ is an isomorphism and that $N_\p \cong (N_\p)^{\dag_\p \dag_\p}$ is totally $C_\p$-reflexive.
Therefore we have $\depth N_\p=\depth R_\p\geq\min\{ n, \dim R_\p\}$.
Hence $N$ satisfies $(S_{n-1})$.
\qed
\end{pf1}

This Theorem gives following corollary.

\begin{cor}
Let $C$ and $C'$ be semidualizing $R$-modules.
Assume that $R$ satisfies the conditions $(S_n)$, $(G_{n-1}^C)$ and $(G_{n-1}^{C'})$.
Then the following conditions are equivalent for an $R$-module $M$.
\begin{enumerate}[\rm (1)]
\item $M$ is $n$-$C$-torsionfree,
\item $M$ is $n$-$C'$-torsionfree,
\item $M$ is $n$-$C$-syzygy,
\item $M$ is $n$-$C'$-syzygy.
\end{enumerate}
\end{cor}

\section{Applications}

In this section, we give some applications of Theorem \ref{tf}.
The first application is to characterize Cohen-Macaulay rings $R$ such that $R_\p$ is Gorenstein for all prime ideals $\p$ of height less than $n$. 
We attempt to solve Tachikawa conjecture and some homological conjectures as a second application.

\subsection{Locally Gorensteinness}

Theorem \ref{lg} gives a characterization Cohen-Macaulay rings that is locally Gorenstein in terms of the Serre's condition.

\begin{thm}\label{lg}
Let $R$ be a Cohen-Macaulay ring with a canonical module $\omega$.
Then the following conditions are equivalent:
\begin{enumerate}[\rm (1)]
\item $R$ satisfies $(G_{n-1}^C)$,
\item $S_n(R)=\Omega^n_C(\mod R)$,
\item $\omega\in\Omega^n_C(\mod R)$.
\end{enumerate}
\end{thm}

\begin{pf}
If $n=0$, then three statements hold.
We assume that $n>0$. 

$(1)\Rightarrow (2)$ 
It is obvious by Theorem \ref {tf}. 

$(2)\Rightarrow (3)$
A canonical module $\omega$ satisfies $(S_n)$, so we have $\omega \in \Omega_C^n(\mod R)$.

$(3)\Rightarrow (1)$
There is an exact sequence
$$
0 \to \omega \to P_C^1 \to P_C^2 \to \cdots \to P_C^n \to M \to 0
$$
such that each $P_C^i$ is a direct summand of direct sum of finite copy of $C$.
For any prime ideal $\p$ of height less than $n$, $(\Omega_C^{n-1}M)_\p$ is a maximal Cohen-Macaulay $R_\p$-module. 
Then the exact sequence $0 \to \omega_\p \to (P_C^1)_\p \to (\Omega_C^{n-1}M)_\p \to 0$ splits.
This indicates $\omega_\p \cong C_\p$.
Thus we have $\id_{R_\p}\, C_\p=\id_{R_\p}\, \omega_\p<\infty$. \qed
\end{pf}

Theorem \ref{lg} recovers a result of Leuschke and Weigand \cite[Lemma 1.4]{LW} by taking as $C=R$ and $n=d$.

\begin{cor}
Let $R$ be a $d$-dimensional Cohen-Macaulay local ring with a canonical module $\omega$.
Then the following conditions are equivalent:
\begin{enumerate}[\rm (1)]
\item $R$ satisfies the $(G_{d-1}^R)$,
\item any maximal Cohen-Macaulay modules are $d$-th syzygy,
\item $\omega$ is $d$-th syzygy.
\end{enumerate}
\end{cor}

\subsection{Some homological conjectures}

In this subsection, we give a partial answer of Tachikawa conjecture and Auslander-Reiten conjecture as an another application of Theorem \ref{tf}.

\begin{thm}\label{main}
We set $d=\sup\{ \ \depth R_\p \ |\  \p \in \Spec R \ \}$, and assume that $0<d<\infty$.
Moreover, we assume that $R$ satisfies the conditions $(S_d)$ and $(ES^C_1)$.
If $\Ext^{>0}(M,C)=0$, then $M$ is totally $C$-reflexive.
\end{thm}

\begin{pf}
We show that $M_\p$ is totally $C_\p$-reflexive by induction on $\height \p$.
Since $R$ satisfies $(G_0^C)$, the case where $\height \p=0$ is obvious.
Next, we consider the case where $\height \p>0$.
By localizing at $\p$, we may assume followings:

$R$ is a noetherian local ring of $\depth d>0$. 
$R$ satisfies $(S_d)$ and $(ES^C_1)$. 
$M$ is an $R$-module such that $\Ext^{>0}_R(M,C)=0$ and that $M_\p$ is totally $C_\p$ reflexive for all nonmaximal prime ideal $\p$.

By replacing $M$ by $\Omega M$ if we need, we may assume $\depth M=d$.
Since $R$ satisfies $(S_d)$ and since $M_\p$ is totally $C_\p$ reflexive for all nonmaximal ideal $\p$, we have $\depth M_\p=\depth R_\p \geq \min\{ d, \dim R_\p\}$.
Therefore $M$ satisfies $(S_d)$.
By Theorem \ref{tf}, $M$ has an $n$-$C$-universal pushforward $0 \to M \to F_C^1 \to \cdots \to F_C^d$.
Putting $M_0=M$ and divide this exact sequence to $0 \to M_{i-1} \to F_C^i \to M_i \to 0$ for each $i=1,2,\dots,d$.
Then one can check that $\Ext^{>0}_R(M_i,C)=0$ for all $i=0,1,\dots,d$.

Suppose that there exists $i$ such that $\depth M_i<d$.
Take such $i$ minimal, we can assume that $M_{i-1}$ satisfies $(S_d)$.
By replacing $M$ by $M_{i-1}$, we may assume $\depth M_1=d-1$.
In this case, we have $\depth M_d=d-d=0$ by depth lemma and this contradicts to $(ES^C_1)$.
Hence we get that $\depth M_i \geq d$ and that $M_{i}$ satisfies $(S_d)$ for every $i$.

We apply above argument repeatedly, we get the following exact sequence:

\begin{center}
\begin{picture}(350,40)

\put(0,22){$0$}
\put(11,26){\vector(1,0){24}}
\put (36,22){$M$}
\put (50,26){\vector(1,0){24}}
\put (75,22){$F_C^1$}
\put (89,26){\vector(1,0){24}}
\put (114,22){$F_C^2$}
\put (128,26){\vector(1,0){24}}
\put (153,22){$\cdots$}
\put (167,26){\vector(1,0){24}}
\put (192,22){$F_C^d$}
\put (206,26){\vector(1,0){36}}
\put (243,22){$F_C^{d+1}$}
\put (265,26){\vector(1,0){24}}
\put (290,22){$F_C^{d+2}$}
\put (313,26){\vector(1,0){24}}
\put (339,23){$\cdots$}

\put (207,21){\vector(1,-1){12}}
\put (230,9){\vector(1,1){12}}
\put (217,0){$M_d$}

\end{picture}
\end{center}
whose $C$-dual sequence is also exact.
This implies that $M$ is totally $C$-reflexive.

\qed
\end{pf}

\begin{cor}
Let $R$ be a generically Gorenstein local ring of depth $d$ and let $M$ be an $R$-module.
Assume that $R$ satisfies $(S_d)$ and $(ES^R_1)$.
If $\Ext^{>0}_R(M,R)=0$, then $M$ is totally reflexive.
\end{cor}

The following corollary is a partial result of the Tachikawa conjecture and more general case has been proved by Avramov, Buchweitz and \c{S}ega \cite{ABS}.

\begin{cor}\cite[Theorem 2.1]{ABS}
Let $R$ be a generically Gorenstein Cohen-Macaulay local ring with canonical module $\omega$.
If $R$ satisfies $(ES^R_1)$ and if $\Ext^{>0}_R(\omega,R)=0$, then $R$ is Gorenstein.
\end{cor}

The notion of almost Gorenstein has been introduced by Barucci and Fr\"{o}berg \cite{BF} and generalized by Goto, Matsuoka and Phuong \cite{GMP} and Goto, Takahashi and Taniguchi \cite{GTT}.
One can easily check that almost Gorenstein local ring is generically Gorenstein. 
Furthermore, Goto, Takahashi and Taniguchi show that any almost Gorenstein Cohen-Macaulay local ring that is not Gorenstein is G-regular \cite[Corollary 4.5]{GTT}.
By using these results, we get the following corollary that is a partial result of the Auslander-Reiten conjecture.

\begin{cor}
Let $R$ be an almost Gorenstein Cohen-Macaulay local ring that is not Gorenstein and let $M$ be an $R$-module.
Assume that $R$ satisfies $(ES^R_1)$.
If $\Ext^{>0}_R(M,R)=0$, then $M$ is free.
\end{cor}

\section{Examples of semidualizing modules}

In this section, we give examples of non-trivial semidualizing modules.
Jorgensen, Leuschke and Sather-Wagstaff \cite{JLS} have been determined the structure of rings which admits non-trivial semidualizing modules. 

\begin{thm}\cite[Theorem 1.1, Theorem 3.1]{JLS}\label{jls}
Let $R$ be a local Cohen-Macaulay ring with a canonical module. Then $R$ admits a semidualizing module that is neither canonical module nor $R$ if and only if there exist a
Gorenstein local ring $S$ and ideals $I_1, I_2 \subset S$ satisfying the following conditions:
\begin{enumerate}[\rm (1)]
\item There is a ring isomorphism $R\cong S/(I_1 + I_2)$,
\item For $j = 1, 2$ the quotient ring $S/I_j$ is Cohen-Macaulay and not Gorenstein,
\item For all $i\in \Z$, we have the following vanishing of Tate cohomology modules:
$$\widehat{\Tor}^S_i (S/I_1, S/I_2) = 0 =\widehat{\Ext}^i_S(S/I_1, S/I_2),$$
\item There exists an integer $c$ such that $\Ext^c_S(S/I_1, S/I_2)$ is not cyclic,
\item For all $i \geq 1$, we have $\Tor^S_i (S/I_1, S/I_2) = 0$; in particular, there is an equality $I_1 \cap I_2 = I_1 I_2$.
\end{enumerate}
Then $C=\Ext^{\dim S-\dim R}_S(R,S)$ is a semidualizing module that is neither $R$ nor canonical module.
\end{thm}

We give a class of Cohen-Macaulay local rings $R$ which have a non-trivial semidualizing module $C$ by using this result.
Moreover, $C_\p$ is a canonical $R_\p$-module for all non-maximal prime ideal $\p$ of $R$.

\begin{prop}
Let $k$ be a field and $T=k[[x_1,x_2,\dots,x_m,y_1,y_2]]$ be a formal power series ring.
For $f_1,f_2,\dots,f_r\in k[[x_1,x_2,\dots,x_m]]$ and $\ell\geq 2$,
we set ideals $I_1=(f_1,f_2,\dots,f_r)T$ and $I_2=(y_1, y_2)^\ell T$.
Assume that $S=T/I_1$ is a $(d+2)$-dimensional Cohen-Macaulay ring which is not Gorenstein and that $S$ satisfies $(G_{n+2}^S)$.
Putting $R=S/I_2$ and $C=\Ext^2_S(R,S)$, then the followings hold:
\begin{enumerate}[\rm (1)]
\item $R$ is $d$-dimensional Cohen-Macaulay ring,
\item $C$ is neither $R$ nor canonical $R$-module,
\item $R$ satisfies $(G_{n}^C)$.
\end{enumerate}
\end{prop}

\begin{pf}
(1) is clear. (2) is comes from Theorem \ref{jls}.
We show (3).
Let $\p$ be a prime ideal of $R$ with height at most $n$.
Since $P=\p S$ is a prime ideal of $S$ with height at most $n+2$, we have that $S_\p=S_P$ is Gorenstein.
Therefore $C_\p=\Ext^2_{S_P}(R_\p,S_P)$ is a canonical $R_\p$-module. \qed
\end{pf}

In the end of this article, we give examples of $d$-dimensional Cohen-Macaulay rings $R$ and semidualizing module $C$ such that $R$ satisfies $(G_{d-1}^C)$ but not $(G_{n}^R)$ for all $n$.

\begin{ex}
Let $k$ be a field.

(1) Let $S=k[[x_1,x_2,x_3,y_1,y_2]]/(x_2^2-x_1x_3, x_2x_3, x_3^2)$ be a $3$-dimensional Cohen-Macaulay local ring which is not Gorenstein.
We set $R=S/(y_1^2, y_1y_2, y_2^2)$ which is a $1$-dimensional Cohen-Macaulay local ring. 
Note that all the prime ideals of $R$ are $\p=(x_2,x_3,y_1,y_2)$ and $\m=(x_1,x_2,x_3,y_1,y_2)$.
It is easy to see that $S_\p$ is Gorenstein but $R_\p$ is not Gorenstein. 
In particular, $R$ does not satisfy $(G_0^R)$.
Putting $C=\Ext_{S}^2(R,S)$, one can check that $C$ is a semidualizing $R$-module which is neither $R$ nor canonical module.
Since $S_\p$ is Gorenstein, we can see that $C_\p$ is a canonical module over $R_\p$.
This yield that $R$ satisfies $(G_0^C)$.

(2) Let $S=k[[x_1,x_2,x_3,x_4,y_1,y_2]]/(x_2^2-x_1x_3, x_2x_3-x_1x_4,x_3^2-x_2x_4)$ be a $4$-dimensional Cohen-Macaulay local normal domain which is not Gorenstein. 
We set $R=S/(y_1^2, y_1y_2, y_2^2)$ which is a $2$-dimensional Cohen-Macaulay local ring. 
For each prime ideal $\p$ except for $\m=(x_1,x_2,x_3,x_4,y_1,y_2)$, we easily show that $S_\p$ is regular. 
But for the unique minimal prime ideal $\q = (x_2^2-x_1x_3, x_2x_3-x_1x_4,x_3^2-x_2x_4,y_1,y_2)$, $R_\q$ is isomorphic to $K[[y_1,y_2]]_{(y_1,y_2)}/(y_1^2,y_1y_2,y_2^2)$, where $K$ is a field. 
Thus $R$ does not satisfy $(G_0^R)$ and therefore $R$ does not satisfy $(G_n^R)$ for every $n$.
Putting $C=\Ext_{S}^2(R,S)$, one can check that $C$ is a semidualizing $R$-module which is neither $R$ nor canonical module. 
Since $S_\p$ is regular, we can see that $C_\p$ is a canonical module over $R_\p$.
This yield that $R$ satisfies $(G_1^C)$.

(3) Let $S=k[[x_1,x_2,x_3,x_4,x_5,x_6,y_1,y_2]]/(x_2^2-x_1x_4, x_2x_3-x_1x_5,x_3^2-x_1x_6, x_3x_4-x_2x_5, x_3x_5-x_2x_6, x_5^2-x_4x_6)$ be a $5$-dimensional Cohen-Macaulay local domain with an isolated singularity which is not Gorenstein. 
We set $R=S/(y_1^2, y_1y_2, y_2^2)$ which is a $3$-dimensional Cohen-Macaulay local ring. 
For each prime ideal $\p$ except for $\m=(x_1,x_2,x_3,x_4,x_5,x_6,y_1,y_2)$, we easily show that $S_\p$ is regular. 
But for the unique minimal prime ideal $\q = (x_2^2-x_1x_4, x_2x_3-x_1x_5,x_3^2-x_1x_6, x_3x_4-x_2x_5, x_3x_5-x_2x_6, x_5^2-x_4x_6,y_1,y_2)$, $R_\q$ is isomorphic to $K[[y_1,y_2]]_{(y_1,y_2)}/(y_1^2,y_1y_2,y_2^2)$, where $K$ is a field. 
Thus $R$ does not satisfy $(G_0^R)$ and therefore $R$ does not satisfy $(G_n^R)$ for every $n$.
Putting $C=\Ext_{S}^2(R,S)$, one can check that $C$ is a semidualizing $R$-module which is neither $R$ nor canonical module. 
Since $S_\p$ is regular, we can see that $C_\p$ is a canonical module over $R_\p$.
This yield that $R$ satisfies $(G_2^C)$.

\end{ex}


\end{document}